# OPTIMAL DESIGNS FOR MIXED MODELS IN EXPERIMENTS BASED ON ORDERED UNITS

By Dibyen Majumdar[1] and John Stufken[2]

*University of Illinois at Chicago and University of Georgia*

We consider experiments for comparing treatments using units that are ordered linearly over time or space within blocks. In addition to the block effect, we assume that a trend effect influences the response. The latter is modeled as a smooth component plus a random term that captures departures from the smooth trend. The model is flexible enough to cover a variety of situations; for instance, most of the effects may be either random or fixed. The information matrix for a design will be a function of several variance parameters. While data will shed light on the values of these parameters, at the design stage, they are unlikely to be known, so we suggest a *maximin* approach, in which a minimal information matrix is maximized. We derive *maximin universally optimal designs* and study their robustness. These designs are based on semibalanced arrays. Special cases correspond to results available in the literature.

**1. Introduction.** When planning an experiment to compare different treatments, it is important that we carefully consider the possible presence of systemic natural differences between the experimental units to be used. If such differences are thought to exist, blocking and the use of covariates are two methods that may help to increase the sensitivity of the experiment for detecting possible differences between the treatments. These two methods are at the core of this paper.

Blocking always leads to a restricted randomization, in which, for each block, a selected set of treatments is randomly assigned to the experimental units in that block. The use of covariates only leads to a restriction on the randomization if the covariates are already used at the design stage rather

Received June 2007; revised June 2007.
[1]Supported by NSF Grant DMS-02-04532.
[2]Supported by NSF Grant DMS-07-06917.

than just at the analysis stage. If covariates are used at the design stage, the designs are often referred to as *systematic designs*, even though there is usually still some opportunity for a restricted randomization.

Cox [5] considered an experiment involving the processing of wool, using a different treatment each week. The natural aging of the wool (which formed the experimental units) caused a trend in the units over time and made time a convenient and useful covariate to account for systemic changes in the experimental units. Cox showed that a systematic assignment of treatments to the units that allowed for estimating the treatment differences in the same way as without the trend was preferable to a fully randomized assignment of the treatments to the units, or to attempting to reduce the effect of the trend by blocking the units.

The covariates that we will consider in this paper are precisely of this type, that is, they are based on a natural ordering of the experimental units, typically induced by time or spatial location. However, our discussion is entirely in the context of block designs, in which each block has the same number $k$ of experimental units. In each block, the units are labeled from 1 through $k$, which induces the covariate (or possibly covariates). The designs in this paper are relevant if it is thought that units may change gradually across this ordering in each of the blocks, although this change may differ from one block to the next. This change could, for example, occur as the result of a learning process, equipment or product deterioration, or spatial location. For some examples and further discussion and references, we refer the reader to Bradley and Yeh [1], Chai and Majumdar [3], Lin and Dean [11], Lin and Stufken [12], Jacroux, Majumdar and Shah [7, 8] and Majumdar and Martin [15]. Lin and Stufken [12] also contains some additional discussion on the pros and cons of using systematic designs.

Thus, we consider the situation where experimental units are partitioned into equally large groups of relatively homogeneous units, or *blocks*, and where, within each block, the units are linearly ordered over time or space. We will build a model that is more flexible than models thus far considered for this situation, and that contains other models as special cases. The model will include random block effects (which contains the model with fixed block effects as a special case), random trend effects that may differ from one block to the next (which contains fixed trend effects, whether the same for each block or varying over the blocks, as a special case) and unit-specific random deviations from the trend (motivated by our belief that a smooth trend is often not very realistic). The latter is a feature hitherto not used in this arena.

While this mixed-effects model will be very flexible, it will also typically contain a considerable number of unknown covariance parameters. Data may help to shed some light on these parameters at the analysis stage, but this is of little help at the design stage. The determination of an optimal design for estimation of the treatment differences, which is the objective of this paper, would therefore seem to be a rather intractable problem since the information matrix for the treatment effects will depend on the many unknown covariance parameters. We will address this problem by identifying, for each design, a "minimal" information matrix for the treatment effects. This minimal information matrix, which will be smaller in the Löwner ordering than the actual information matrix for the treatment effects, will depend on very few (no more than two) parameters, which are functions of the original covariance parameters. It is this minimal information matrix that we will *maximize* over all designs to obtain a "maximin" information matrix and an optimal design. We note that our approach is in the spirit of Kiefer [9] and Kiefer and Wynn [10], who considered minimax optimal designs for models with autoregressive errors.

The maximin information matrix will be derived in Section 2. After characterizing and identifying the optimal designs in Section 3 (with proofs deferred to Section 5), we will investigate the robustness of this process to parameter misspecification in Section 4.

**2. Setup and basic results.** Consider an experiment to compare $v$ treatments ($i = 1, \ldots, v$) based on $n = bk$ experimental units that are partitioned into $b$ blocks ($j = 1, \ldots, b$) of $k$ units each ($p = 1, \ldots, k$). Suppose the units within blocks are linearly ordered over time or space. The collection of units can be visualized as a $k \times b$ array with rows labeled by units within blocks and columns by blocks, while the entries of the array are treatments assigned to the units by the design, that is, for design $d$, the entry in cell $(p, j)$ is $d(p, j)$, where $d(p, j) \in \{1, \ldots, v\}$. The design $d$ itself will be viewed as a $k \times b$ matrix. Under a very general design, our objective is to determine an optimal design for comparing the treatments.

For the model, in addition to a *block effect*, we assume that there is a trend over time or space within each block. If $y_{pj}$ denotes the random variable corresponding to the observation in row (unit) $p$ and column (block) $j$, then the model is

$$y_{pj} = \mu + \tau_{d(p,j)} + \beta_j^b + \zeta_{pj} + \varepsilon_{pj},$$

where $\tau_{d(p,j)}$ denotes an effect for treatment $d(p, j)$, $\beta_j^b$ an effect for block $j$, $\zeta_{pj}$ a trend effect for unit $p$ in block $j$ and $\varepsilon_{pj}$ the measurement error. We assume that the trend $\zeta_{pj}$ is composed of two parts, one that is smooth enough to be approximated by a polynomial in $p$ and another that represents random fluctuations from the polynomial. The smooth part, which we assume to be linear, will be written as $\gamma_j^I + \gamma_j \phi(p)$, where $\phi(p)$ is the linear orthonormal polynomial on $\{1, \ldots, k\}$, specifically, $\phi(p) = \sqrt{3/(k(k^2-1))}(2p - k - 1)$. If $\psi_{pj}$ denotes the fluctuation from the smooth trend then we may write $\zeta_{pj} = \gamma_j^I + \gamma_j \phi(p) + \psi_{pj}$. The slope $\gamma_j$ is further decomposed into a fixed part ($\theta_0$) that is common to all blocks and a random part ($\theta_j$) that may vary from block to block, that is, $\gamma_j = \theta_0 + \theta_j$. Writing $\beta_j = \beta_j^b + \gamma_j^I$ and $\delta_{pj} = \psi_{pj} + \varepsilon_{pj}$, we arrive at our model,

$$(2.1) \qquad y_{pj} = \mu + \tau_{d(p,j)} + \beta_j + \theta_0 \phi(p) + \theta_j \phi(p) + \delta_{pj}.$$

This is a mixed-effects model. The quantities $\tau_1, \ldots, \tau_v$ and $\theta_0$ are considered to be fixed effects while $\beta_j$, $\theta_j$, $\delta_{pj}$ for $p = 1, \ldots, k$, $j = 1, \ldots, b$ are considered random. Although (2.1) combines the two variables $\psi_{pj}$ and $\varepsilon_{pj}$ into one variable $\delta_{pj}$, at the modeling stage, it may be useful to recognize their individual characteristics. Generally, the measurement error $\varepsilon_{pj}$ may be assumed to be homoscedastic, while $\psi_{pj}$, the departure from the assumed trend, will depend on the strength, nature and variability of the trend in the particular application. These considerations may enable the experimenter to determine an appropriate variance–covariance structure for the $\delta_{pj}$'s.

We assume that the random effects have zero expectations and are uncorrelated from one block to the next. Let $\sigma_\beta^2$ and $\sigma_\theta^2$ denote the variances of $\beta_j$ and $\theta_j$, $\sigma_{\beta\theta}$ their covariance, $V_{\delta\beta}$ and $V_{\delta\theta}$ the $k \times 1$ vectors of covariances of $\delta_j = (\delta_{1j}, \ldots, \delta_{kj})'$ with $\beta_j$ and $\theta_j$ and $V_{\delta\delta}$ the covariance matrix of the $\delta$'s. Let $\mathbf{1}_k$ denote the $k \times 1$ vector of 1's, $\tau = (\tau_1, \ldots, \tau_v)'$ and $\phi = (\phi(1), \ldots, \phi(k))'$. If $Y_j$ denotes the $k \times 1$ vector of observations from block $j$, then $E(Y_j) = X_{dj}\tau + Z_0\gamma$, where $\gamma = (\mu, \theta_0)'$, $Z_0 = (\mathbf{1}_k, \phi)$ and $X_{dj}$ is the $k \times v$ unit-treatment incidence matrix for block $j$ with entries

$$X_{dj}(p, i) = \begin{cases} 1, & \text{if } d(p,j) = i, \\ 0, & \text{otherwise.} \end{cases}$$

Also, $V(Y_j) = \Sigma = \sigma_\beta^2 \mathbf{1}_k \mathbf{1}_k' + \sigma_\theta^2 \phi\phi' + V_{\delta\delta} + \sigma_{\beta\theta}(\mathbf{1}_k \phi' + \phi \mathbf{1}_k') + (\mathbf{1}_k V_{\delta\beta}' + V_{\delta\beta} \mathbf{1}_k') + (\phi V_{\delta\theta}' + V_{\delta\theta} \phi')$. Let $Z = (\mathbf{1}_{bk}, \mathbf{1}_b \otimes \phi)$ and $X_d = (X_{d1}', \ldots, X_{db}')'$, where $\otimes$ denotes Kronecker product. The first and second moments of the observations $Y = (Y_1', \ldots, Y_b')'$ are then

$$(2.2) \qquad E(Y) = X_d\tau + Z\gamma, V(Y) = V = \mathbf{I}_b \otimes \Sigma,$$

where $\mathbf{I}_b$ is the identity matrix of order $b$. The parameter of interest is $\tau$, the vector of treatment effects. For a design $d$, the information matrix for $\tau$ is given by

$$C_d = X_d' V^{-1} X_d - X_d' V^{-1} Z (Z' V^{-1} Z)^{-1} Z' V^{-1} X_d.$$

REMARK 2.1 [Special cases of (2.2)]. If $V_{\delta\delta} = a_0 \mathbf{I}_k$ with $a_0 > 0$ (equivalently, $V_{\delta\delta} = a_0 \mathbf{I}_k + a_1 \mathbf{1}_k \mathbf{1}_k' + a_2 \phi\phi'$ with $a_0 > 0$), then the model is equivalent to one in which the random trend in each block is known to be linear. If, in addition, $\sigma_\theta^2 = 0$, then the model is equivalent to one in which the trend is fixed, linear and the same for each block. This has been studied by several authors, including Bradley and Yeh [1], Yeh and Bradley [21], Yeh, Bradley and Notz [22], Stufken [20], Lin and Dean [11], Chai and Majumdar [3], Chai [2] and Lin and Stufken [12, 13]. Alternatively, if $\sigma_\theta^2 = \infty$, then the model is still equivalent to one in which the trend is fixed and linear, but now possibly different in different blocks. This has been studied by Jacroux, Majumdar and Shah [7, 8] and Majumdar and Martin [14]. For any $\sigma_\theta^2$, if $\sigma_\beta^2 = \infty$, then the model corresponds to one in which the block effects are fixed.

Our objective is to identify a *universally optimal design* for estimating $\tau$. For (2.2), $C_d$ depends on the $4 + 2k + k(k+1)/2$ variance component parameters in $\Sigma$. If these are all known at the planning stage, then the optimization problem can be solved by numerical techniques. However, this will rarely be the case. Our approach, therefore, is to work with few parameters at the design stage. To do this, we first identify, for each design $d$, a minimal information matrix $C_d^L$, and then determine a universally optimal design based on the minimal information matrix. This is, therefore, a maximin approach. We will derive a minimal information matrix that depends on at most two parameters (other than $v$, $b$ and $k$), which are functions of the original variance components. Once the data has been collected, however, we recommend a less parsimonious approach. At the inference stage, the experimenter should work with a realistic model with all likely parameters included and let the data decide.

We will use the Löwner ordering to identify the minimal information matrix, that is, $B \succeq A$ if $B - A$ is nonnegative definite. Formally, the first step is to identify, for each $d$, a matrix $C_d^L$ such that $C_d \succeq C_d^L$, where $C_d^L$ is an information matrix for the design $d$ corresponding to a simplified model. The next step is to find the optimal design $d^*$ such that

$$C_{d^*}^L \text{ is completely symmetric} \quad \text{and} \quad \text{trace}(C_{d^*}^L) = \max_{d \in D}(C_d^L).$$

To get $C_d^L$, we utilize the representation $C_d = X_d' Q(QVQ)^- Q X_d$, where $Q = \mathbf{I}_n - Z(Z'Z)^{-1}Z'$, and observe that a lower bound for $C_d$ can be obtained by using an upper bound for $\Sigma$. To get an upper bound for $\Sigma$, we note that in most situations, variances are easier to determine than covariances. It can be shown that $\Sigma = V(\delta_j + \beta_j \mathbf{1}_k + \theta_j \phi) \preceq 4V(\delta_j) + 4V(\beta_j \mathbf{1}_k) + 4V(\theta_j \phi) \preceq 4E_{\max}(V_{\delta\delta})\mathbf{I}_k + 4\sigma_\beta^2 \mathbf{1}_k \mathbf{1}_k' + 4\sigma_\theta^2 \phi\phi'$, where $E_{\max}(V_{\delta\delta})$ is the maximum eigenvalue of $V_{\delta\delta}$. This bound is generally conservative, and we can do better if there is additional information. For example, if $V_{\delta\beta} = V_{\delta\theta} = 0$, that is, the fluctuations from the trend $\delta_j$ are not correlated with the modeled part of the trend $\theta_j$ and the block effect $\beta_j$, then $\Sigma \preceq E_{\max}(V_{\delta\delta})\mathbf{I}_k + 2\sigma_\beta^2 \mathbf{1}_k \mathbf{1}_k' + 2\sigma_\theta^2 \phi\phi'$.
Hence, in general, we have

$$\Sigma \preceq \sigma_{0\varepsilon}^2 \mathbf{I}_k + \sigma_{0\beta}^2 \mathbf{1}_k \mathbf{1}_k' + \sigma_{0\theta}^2 \phi\phi',$$

using quantities $\sigma_{0\varepsilon}^2, \sigma_{0\beta}^2$ and $\sigma_{0\theta}^2$ that take values in the intervals

$$E_{\max a}(V_{\delta\delta}) \leq \sigma_{0\varepsilon}^2 \leq 4E_{\max a}(V_{\delta\delta}), \qquad \sigma_{\beta a}^2 \leq \sigma_{0\beta}^2 \leq 4\sigma_{\beta a}^2,$$

$$\sigma_{\theta a}^2 \leq \sigma_{0\theta}^2 \leq 4\sigma_{\theta a}^2,$$

where $\sigma_{\beta a}^2$, $\sigma_{\theta a}^2$ and $E_{\max a}(V_{\delta\delta})$ are the assumed or prior values of $\sigma_\beta^2$, $\sigma_\theta^2$ and $E_{\max}(V_{\delta\delta})$, respectively. If the correlations are believed to be negligible the values of $\sigma_{0\varepsilon}^2, \sigma_{0\beta}^2$ and $\sigma_{0\theta}^2$ should be taken at the lower endpoints of the intervals or close to it, while for stronger correlations, these values should be assumed higher. We will see in Section 4 that the optimal designs are remarkably robust, so an accurate determination of $\sigma_{0\varepsilon}^2, \sigma_{0\beta}^2$ and $\sigma_{0\theta}^2$ within their respective intervals is usually not necessary.

Our first theorem gives the minimal information matrix. To state it, we use the standard notation for a design $d$: $r_{di}$ will denote the replication of treatment $i$, $\mathbf{r}_d = (r_{d1}, \ldots, r_{dv})'$, $R_d = \mathrm{diag}(r_{d1}, \ldots, r_{dv})$, $N_d = (n_{dij})$ the treatment $\times$ block (column) incidence matrix and $M_d = (m_{dip})$ the treatment $\times$ unit (row) incidence matrix. Also, let

$$(2.3) \qquad \lambda_0 = \frac{\sigma_{0\beta}^2}{\sigma_{0\varepsilon}^2 + k\sigma_{0\beta}^2}, \qquad \lambda_1 = \frac{\sigma_{0\theta}^2}{\sigma_{0\varepsilon}^2 + \sigma_{0\theta}^2}.$$

THEOREM 2.2. *The information matrix $C_d$ for a design $d$ based on the model (2.2) satisfies $C_d \succeq C_d^L$, where*

$$(2.4) \quad \sigma_{0\varepsilon}^2 C_d^L = \sum_{j=1}^b X_{dj}' W X_{dj} - \left(\frac{1-k\lambda_0}{bk}\right)\mathbf{r}_d \mathbf{r}_d' - \left(\frac{1-\lambda_1}{b}\right) M_d \phi\phi' M_d'$$

*with*

$$(2.5) \qquad W = \mathbf{I}_k - \lambda_0 \mathbf{1}_k \mathbf{1}_k' - \lambda_1 \phi\phi'.$$

The proof is straightforward and hence omitted.

REMARK 2.3. (i) $C_d^L$ is the information matrix for model (2.2) with uncorrelated random effects ($V_{\delta\beta} = V_{\delta\theta} = 0, \sigma_{\beta\theta} = 0$), $V(\delta_j) = V_{\delta\delta} = \sigma_{0\varepsilon}^2 \mathbf{I}_n$, $V(\beta_j) = \sigma_\beta^2 = \sigma_{0\beta}^2$ and $V(\theta_j) = \sigma_\theta^2 = \sigma_{0\theta}^2$.

(ii) The minimal information matrix may also be written as

$$
\sigma_{0\varepsilon}^2 C_d^L = R_d - \lambda_0 N_d N_d' - \lambda_1 \sum_{j=1}^{b} X_{dj}' \phi \phi' X_{dj}
$$
$$
- \left(\frac{1-k\lambda_0}{bk}\right) \mathbf{r}_d \mathbf{r}_d' - \left(\frac{1-\lambda_1}{b}\right) M_d \phi \phi' M_d'.
$$

**3. Optimal designs.** In this section, we will explore optimal designs for model (2.1). Our goal is to determine universally optimal designs using the minimal information matrix (2.4), that is, a maximin universally optimal design. First, we need a definition.

DEFINITION 3.1 (Rao [18, 19]). A $t \times b$ array in $v$ symbols is called an *orthogonal array of type II (OAII) of strength* 2 or a *semibalanced array* if the columns of the array consist of distinct symbols and any two rows of the array contain every pair of distinct symbols equally often.

For the construction of these arrays, see Hedayat, Sloane and Stufken [6].

We will establish the maximin universal optimality of judiciously selected designs of the form

$$
(3.1) \qquad \widetilde{d} = \Pi \left( \frac{\widetilde{d_1}}{\widetilde{d_2}} \right),
$$

where $\widetilde{d_1}$ is a semibalanced array with rows that are uniform in symbols, $\widetilde{d_2}$ is a matrix with each row identical to some row of $\widetilde{d_1}$ and $\Pi$ is a permutation matrix of order $k$. Note that the rows of a semibalanced array are always uniform when there are three or more rows.

Universality optimality will be established by the technique outlined in Theorem 1 of Majumdar and Martin [14]. To establish complete symmetry (c.s.), note that each treatment occurs equally often in each row of $\widetilde{d}$. Hence, for a design of the form (3.1), $M_{\widetilde{d}}\phi = 0$ and $\mathbf{r}_{\widetilde{d}}\mathbf{r}_{\widetilde{d}}'$ is c.s. Also, it follows from Cheng [4] that $\sum_{j=1}^{b} X_{\widetilde{d}j}' W X_{\widetilde{d}j}$ is c.s. Hence, $C_{\widetilde{d}}^L$ in (2.4) is c.s.

For an arbitrary design $d \in D$, the trace of the maximin information matrix is given by

$$
\sigma_{0\varepsilon}^2 \operatorname{trace}(C_d^L) = \operatorname{trace}\left( \sum_{j=1}^{b} X_{dj}' W X_{dj} \right) - \left(\frac{1-k\lambda_0}{bk}\right) \mathbf{r}_d' \mathbf{r}_d
$$
$$
- \left(\frac{1-\lambda_1}{b}\right) \phi' M_d' M_d \phi.
$$

Note that $\phi' M_d' M_d \phi \geq 0$, $r_d' r_d \geq \frac{(bk)^2}{v}$ and both of these lower bounds are attained by designs of the form (3.1). Since

$$(3.2) \qquad\qquad 0 \leq \lambda_0 \leq \frac{1}{k}, \qquad 0 \leq \lambda_1 \leq 1,$$

if there is a design of the form (3.1) that maximizes trace($\sum_{j=1}^{b} X_{dj}' W X_{dj}$) among all designs, then this is maximin universally optimal.

Writing, $W = (w_{pq})$ we get trace($X_{dj}' W X_{dj}$) $= \sum_{p=1}^{k} w_{pp} + 2 \sum_{(p,q)} w_{pq}$, where the second summation is over all pairs of experimental units $p, q$, $p < q$, that are occupied by the same treatment in block $j$. Since $\sum_{p=1}^{k} w_{pp}$ does not depend on the design, an order of assignment of treatments to block $j$ that maximizes the second sum will maximize trace($X_{dj}' W X_{dj}$). Note that since $w_{pq}$ does not depend on the block subscript $j$, the pattern that is optimal for block $j$ is also optimal for any other block and, combined, will maximize trace($\sum_{j=1}^{b} X_{dj}' W X_{dj}$).

Let $\mathcal{O} = \{\pi = (\pi(1), \ldots, \pi(k)) : \pi(p) \in \{1, \ldots, v\}, \ p = 1, \ldots, k\}$ be the set of orders $\pi$. For each order $\pi \in \mathcal{O}$, there is a design of the form (3.1) in which $\pi$ is the first block as long as a $k^* \times b$ semibalanced array based on $v$ treatments exists, where $k^*$ is the number of distinct treatments in $\pi$. For $\pi \in \mathcal{O}$, let

$$P(\pi) = \{(p, q) : 1 \leq p < q \leq k, \pi(p) = \pi(q)\},$$
$$F(\pi) = \sum_{p, q \in P(\pi)} w_{pq}.$$

Note that $F(\pi)$ depends on $v, k, \lambda_0$ and $\lambda_1$, but not on $b$.

The trace maximization problem may be stated as, given $v, k, \lambda_0$ and $\lambda_1$ satisfying (3.2), maximize $F(\pi)$ over all $\pi \in \mathcal{O}$. An order $\pi^*$ that maximizes $F(\pi)$ will be called an *optimal order*. In the next two subsections, we will identify optimal orders $\pi^*$. We will distinguish between the two cases $k < 2v$ and $k \geq 2v$. The proofs are given in Section 5.

Before concluding this section, we give an alternate expression for $F(\pi)$. Since for $p \neq q$, $w_{pq} = -\lambda_0 - \lambda_1 \phi(p) \phi(q)$, we have

$$(3.3) \qquad\qquad F(\pi) = -\lambda_0 s(\pi) - \lambda_1 T(\pi),$$

where $s(\pi) = |P(\pi)|$, the cardinality of $P(\pi)$, and $T(\pi) = \sum_{p, q \in P(\pi)} \phi(p) \phi(q)$. For $i = 1, \ldots, v$, if we denote

$$(3.4) \qquad\qquad n_i = n_i(\pi) = \text{ replication of treatment } i \text{ in } \pi,$$

$$(3.5) \qquad\qquad h_i = h_i(\pi) = \sum_{p=1}^{k} \delta_{ip} \phi(p),$$

where $\delta_{ip} = \delta_{ip}(\pi) = 1$ if $\pi(p) = i$ and equals zero otherwise, then it follows from (3.4) that

$$(3.6) \qquad s = s(\pi) = \sum_{i=1}^{v} \frac{n_i(n_i - 1)}{2} = \frac{1}{2}\left[\sum_{i=1}^{v} n_i^2 - k\right].$$

Also, since

$$2T(\pi) = \sum_{p,q \in P(\pi)} 2\phi(p)\phi(q) = \sum_{i=1}^{v}\left(\sum_{p=1}^{k} \delta_{ip}\phi(p)\right)^2 - \sum_{p=1}^{k} \phi^2(p),$$

using (3.5), we get

$$(3.7) \qquad T(\pi) = \frac{1}{2}\left[\sum_{i=1}^{v} h_i^2 - 1\right].$$

3.1. *Optimal orders when $k < 2v$.* For a positive integer $q \geq k - v$, we define $\pi_q$ to be an order of the form

$$(3.8) \qquad \pi_q = \{i_1, i_2, \ldots, i_q, i_{q+1}, \ldots, i_{k-q}, i_q, \ldots, i_2, i_1\},$$

where $i_1, i_2, \ldots, i_{k-q}$ are $k - q$ distinct treatments. We will also write $\pi_0$ for an order of $k$ distinct treatments. We define

$$(3.9) \qquad s^* = \begin{cases} \text{Max}\left\{p : p \text{ integer}, 1 \leq p < \dfrac{k+1}{2}, \lambda_1\phi^2(p) > \lambda_0\right\}, \\ \qquad \text{if } \lambda_1\phi^2(1) > \lambda_0, \\ 0, \qquad \text{if } \lambda_1\phi^2(1) \leq \lambda_0. \end{cases}$$

Note that $s^*$ is a function of $k, \lambda_0$ and $\lambda_1$, but not of $v$ or $b$, and $s^* \leq k/2$.

LEMMA 3.2. *Suppose $k < 2v$.*

(i) *If $k \leq v + s^*$, then $\pi_{s^*}$ is an optimal order.*
(ii) *If $k > v + s^*$, let $\alpha = k - v$. $\pi_\alpha$ is then an optimal order.*

Using (3.6) and (3.7), we note that $\pi_\alpha$ minimizes $s(\pi)$ over $\pi \in \mathcal{O}$ and minimizes $T(\pi)$ among all orders that minimize $s(\pi)$. A proof of Lemma 3.2 is given in Section 5.

3.2. *Optimal orders when $k \geq 2v$.* For given $k$ and $v$, $k \geq 2v$, we define integers $m$ and $t$ by

$$(3.10) \qquad k = mv + t \qquad \text{where } 0 \leq t < v \text{ and } m \geq 2.$$

Depending on the values of $m$ and $t$, an optimal order will turn out to be either a *trend-free* (TF) or *nearly trend-free* (NTF) order. An order is called *trend-free* if all treatments are "orthogonal" to the fixed part of the trend, that is, for each $i = 1, \ldots, v$,

$$(3.11) \qquad\qquad h_i = 0.$$

It is easy to see that a trend-free order can only exist if, for each $i = 1, \ldots, v$,

$$(3.12) \qquad\qquad n_i(k+1) \equiv 0 \pmod{2}.$$

When $k$ is odd, any integer $n_i$ satisfies (3.12). A TF order can be constructed in this case if $n_i \geq 2$ for all $i$. The treatments with even replication are used at the beginning and at the end of such an order $\pi$ in such a way that $\pi(p) = \pi(k - p + 1)$ for these positions. For treatments with odd $n_i$, $n_i \geq 3$, (3.11) can be achieved by filling the remaining positions using the construction of Phillips [17], which is also reproduced in Lemma 3.2(a) of Jacroux, Majumdar and Shah [8]. It follows from (3.7) that for a trend-free order, $T(\pi) = -1/2$, the lower bound.

When $k$ is even, (3.12) implies that $n_i$ must be even for all $i$, in which case an order $\pi$ with the property $\pi(p) = \pi(k - p + 1)$ for all $p$ satisfies (3.11). If $n_i$ is odd for some $i$, then it can be shown (see, e.g., Lemma 3.1 of Jacroux, Majumdar and Shah [8]) that

$$(3.13) \qquad\qquad |h_i| \geq \phi\left(\frac{k+2}{2}\right) = -\phi\left(\frac{k}{2}\right).$$

Provided that $n_i \geq 3$ for treatments with an odd replication, using the construction of Mukerjee and Sengupta [16], which is also reproduced in Lemma 3.2(b) of Jacroux, Majumdar and Shah [8], the lower bound in (3.13) can be achieved for each such treatment, while at the same time achieving $h_i = 0$ for treatments with even replications. For $k$ even, orders that satisfy (3.11) and the lower bound in (3.13) for treatments with even $n_i$ and odd $n_i$, respectively, are called *nearly trend-free* orders. Note that for a nearly trend-free order $\pi$ in which $u$ treatments have odd replications, $T(\pi) = \frac{1}{2}[u(\phi(\frac{k}{2}))^2 - 1]$.

The trend-free and nearly trend-free orders that are especially relevant to our investigation are given in the following definitions.

DEFINITION 3.3 (Trend-free orders).

(i) When $k$ is odd, we use $\pi_{TF}^A$ to denote any order with the following two properties.

(a) The replications are $n_1 = \cdots = n_t = m + 1$, $n_{t+1} = \cdots = n_v = m$, there being no treatment with replication $m+1$ if $t = 0$. Treatments with even replication occupy the "outer" positions [which are $1, \ldots, (m+1)t/2$ and $k - (m+1)t/2 + 1, \ldots, k$ when $m$ is odd and $1, \ldots, m(v-t)/2$ and $k - m(v-t)/2 + 1, \ldots, k$ when m is even]. The remaining, "inner," positions are occupied by treatments with odd replication.

(b) Treatments with even replication are arranged so that $\pi_{TF}^A(p) = \pi_{TF}^A(k-p+1)$. Treatments with odd replication are arranged using the construction of Phillips [17] so that $h_i = 0$.

(ii) When $k$ is even and $k/v$ is an even integer (so that $m$ is even and $t = 0$), we use $\pi_{TF}^B$ to denote any order with $n_i = m$, $i = 1, \ldots, v$, and $\pi_{TF}^B(p) = \pi_{TF}^B(k-p+1)$ for $p = 1, \ldots, k/2$.

(iii) When $k$ is even and $k/v$ is not an even integer, we use $\pi_{TF}^C$ to denote any order with $n_i \in \{\xi, \xi+2\}$ for all $i = 1, \ldots, v$, where $\xi$ is the even integer in $\{m-1, m\}$ and $\pi_{TF}^C(p) = \pi_{TF}^C(k-p+1)$ for $p = 1, \ldots, k/2$.

DEFINITION 3.4 (Nearly trend-free orders). When $k$ is even and $k/v$ is not an even integer, we use $\pi_{NTF}$ to denote any order with the following two properties.

(a) The replications are $n_1 = \cdots = n_t = m+1$, $n_{t+1} = \cdots = n_v = m$, there being no treatment with replication $m+1$ if $t = 0$. Treatments with even replication occupy the "outer" positions [which are $1, \ldots, (m+1)t/2$ and $k - (m+1)t/2 + 1, \ldots, k$ when $m$ is odd and $1, \ldots, m(v-t)/2$ and $k - m(v-t)/2 + 1, \ldots, k$ when m is even]. The remaining, "inner," positions are occupied by treatments with odd replication.

(b) Treatments with even replication are arranged so that $\pi_{NTF}(p) = \pi_{NTF}(k-p+1)$. Treatments with odd replication are arranged using the construction of Mukerjee and Sengupta [16] so that $|h_i| = \phi\left(\frac{k+2}{2}\right) = -\phi\left(\frac{k}{2}\right)$.

Optimal orders are given in the following lemma, the proof of which is again postponed to Section 5.

LEMMA 3.5. *Suppose $k \geq 2v$.*

  (i) *If $k$ is odd, then $\pi_{TF}^A$ is an optimal order.*
 (ii) *If $k$ is even and $k/v$ is an even integer, then $\pi_{TF}^B$ is an optimal order.*
(iii) *If $k$ is even and $k/v$ is not an even integer, then*

$$\pi_{TF}^C \text{ is optimal} \qquad \text{if } \lambda_1 \phi^2\left(\frac{k}{2}\right) > \lambda_0;$$

$$\pi_{NTF} \text{ is optimal} \qquad \text{if } \lambda_1 \phi^2\left(\frac{k}{2}\right) \leq \lambda_0.$$

3.3. *Optimal designs.* The main result is formulated in Theorem 3.6 and is an immediate consequence of the considerations in the previous subsections.

THEOREM 3.6. *Given $v$, $k$ and $\lambda_0$, $\lambda_1$ satisfying (3.2), suppose $\pi^*$ is an optimal order given by Lemma 3.2 or 3.5. Suppose $b$ is such that a $k^* \times b$* *semibalanced array based on $v$ treatments exists, where $k^*$ is the number of distinct treatments in $\pi^*$. Let*

$$\widetilde{d}^* = \Pi \left( \frac{\widetilde{d}_1^*}{\widetilde{d}_2^*} \right)$$

*be a $k \times b$ array, where $\widetilde{d}_1^*$ is a $k^* \times b$ semibalanced array, each row of $\widetilde{d}_2^*$ is identical to some row of $\widetilde{d}_1^*$ and $\Pi$ is a permutation matrix such that, after relabeling treatments if necessary, $\widetilde{d}^*(p,1) = \pi^*(p)$, $p = 1, \ldots, k$. $\widetilde{d}^*$ is then a maximin universally optimal design.*

REMARK 3.7. For the case $\sigma_{\beta\theta} = 0$, $V_{\delta\beta} = V_{\delta\theta} = 0$, $V_{\delta\delta} = \sigma_{0\varepsilon}^2 \mathbf{I}_n$, $\sigma_{0\beta}^2 = \infty$ ($\lambda_0 = 1/k$) and $\sigma_{0\theta}^2 = \infty$ ($\lambda_1 = 1$), our results reduce to the results of Jacroux, Majumdar and Shah [8]. Therefore, Theorem 3.6 may be viewed as a generalization of their Corollary 4.3 and Theorem 4.6. In particular, our results extend the results of Jacroux, Majumdar and Shah [13] to models with random block and trend effects. Moreover, our proofs are different from theirs and arguably less cumbersome.

REMARK 3.8. For the case $\sigma_{\beta\theta} = 0$, $V_{\delta\delta} = V_{\delta\theta} = 0$, $V_{\delta\delta} = \sigma_{0\varepsilon}^2 \mathbf{I}_n$, $\sigma_{0\beta}^2 = \infty$ ($\lambda_0 = 1/k$) and $\sigma_{0\theta}^2 = 0$ ($\lambda_1 = 0$), Theorem 3.7 of [3], established the existence and optimality of "strongly balanced" BIB designs. Our approach can be used to generalize this result to models with an arbitrary $\sigma_{0\beta}^2 > 0$ ($\lambda_0 \in (0, 1/k]$), that is, models with random block effects.

## 4. Robustness.

For given $v$, $b$ and $k$, the existence and construction of the optimal designs in Section 3 may require knowledge of the covariance parameters $\lambda_0$ and $\lambda_1$. An important issue is whether the misspecification of these parameters can lead to the choice of inefficient designs. We restrict ourselves to the case $\lambda_1 > 0$. As in Section 3, we will distinguish between the cases $k < 2v$ and $k \geq 2v$, starting with the slightly simpler latter case.

For $k \geq 2v$, if $k$ is odd, or $k$ is even and $k/v$ is an even integer, then the optimal design given by Theorem 3.6 does not depend on $\lambda_0$ or $\lambda_1$. Hence, provided $b$ is such that it accommodates the optimal design in the theorem, for these cases, there is no need to specify $\lambda_0$ or $\lambda_1$ to select an optimal design. On the other hand, if $k$ is even and $k/v$ is not an even integer, then the order for the optimal design in Theorem 3.6 is

$$\pi_{TF}^C \quad \text{if } \frac{\lambda_0}{\lambda_1} < \phi^2 \left( \frac{k}{2} \right) \quad \text{and} \quad \pi_{NTF} \quad \text{if } \frac{\lambda_0}{\lambda_1} \geq \phi^2 \left( \frac{k}{2} \right).$$

Thus, misspecification of $\lambda_0$ or $\lambda_1$ could lead to the selection of $\pi_{TF}^C$ in cases where the design based on $\pi_{NTF}$ is optimal, or vice versa. How bad can this be?

Since $k$ is normally rather large for this case, $\phi^2(\frac{k}{2}) = (2\binom{k+1}{3})^{-1}$ will tend to be small. Therefore, unless $\lambda_0$ is near zero, which corresponds to the case of no block effects, the optimal design will be based on $\pi_{NTF}$. Moreover, that design turns out to be very efficient when it is not optimal. To see this, a natural measure of the relative efficiency of the design based on $\pi_{NTF}$ when the design based on $\pi_{TF}^C$ is optimal is the ratio

$$(4.1) \qquad E_1 = \frac{\text{trace}(C_d^L(\pi_{NTF}))}{\text{trace}(C_d^L(\pi_{TF}^C))},$$

where $C_d^L(\pi)$ stands for the matrix $C_d^L$ for a design as in equation (3.1) that is based on the order $\pi$. When the design based on $\pi_{NTF}$ is optimal, we can use $1/E_1$ to measure the relative efficiency of the design based on $\pi_{TF}^C$.

Expressions for the two traces in (4.1) can easily be computed from the results in Section 3. This leads to

$$\sigma_{0\varepsilon}^2 \, \text{trace}(C_d^L(\pi_{NTF})) = \begin{cases} b\left(k - k\lambda_0 - 2\lambda_0 s^{\min} - \dfrac{k}{v}(1 - k\lambda_0)\right) - bt\lambda_1\phi^2\left(\dfrac{k}{2}\right), \\ \qquad\qquad\qquad\qquad\qquad\qquad\qquad\qquad \text{if } m \text{ is even,} \\ b\left(k - k\lambda_0 - 2\lambda_0 s^{\min} - \dfrac{k}{v}(1 - k\lambda_0)\right) \\ \qquad\qquad\qquad - b(v - t)\lambda_1\phi^2\left(\dfrac{k}{2}\right), \qquad \text{if } m \text{ is odd,} \end{cases}$$

$$\sigma_{0\varepsilon}^2 \, \text{trace}(C_d^L(\pi_{TF}^C)) = \begin{cases} b\left(k - k\lambda_0 - 2\lambda_0 s^{\min} - \dfrac{k}{v}(1 - k\lambda_0)\right) - bt\lambda_0, \\ \qquad\qquad\qquad\qquad\qquad\qquad\qquad\qquad \text{if } m \text{ is even,} \\ b\left(k - k\lambda_0 - 2\lambda_0 s^{\min} - \dfrac{k}{v}(1 - k\lambda_0)\right) - b(v - t)\lambda_0, \\ \qquad\qquad\qquad\qquad\qquad\qquad\qquad\qquad \text{if } m \text{ is odd,} \end{cases}$$

where $m$ and $t$ are defined in (3.10) and

$$(4.2) \qquad s^{\min} = \frac{1}{2}[(v - t)m^2 + t(m + 1)^2 - mv - t] = \frac{m}{2}[k - v + t].$$

Note that $E_1$ does not depend on the value of $b$. It is immediately clear that $E_1$ is virtually equal to 1 if $\lambda_0 = 0$, so the design based on $\pi_{NTF}$ is often optimal and always highly efficient. It is also seen from these expressions that the design based on $\pi_{TF}^C$ is highly efficient when it is not optimal.

Turning to the case $k < 2v$, based on Lemma 3.2, we arrive at an optimal sequence of the form $\pi_q$ in (3.8), where $\max\{0, k - v\} \le q \le \lfloor k/2 \rfloor$. But, if $\lambda_0$ or $\lambda_1$ are misspecified, then we could end up selecting a design with the wrong value of $q$. To see how bad this could be, first note that for all values of $\lambda_0/\lambda_1$ within each of the following intervals, there is one order that is optimal: $(0, \phi^2(\lfloor k/2 \rfloor)]$, $[\phi^2(q + 1), \phi^2(q)]$ for $q = \max\{0, k - v\}, \ldots, \lfloor k/2 \rfloor - 1$, and $[\phi^2(\max\{0, k - v\} + 1), \infty)$. Thus, if the misspecified value and the true

Table 1
*Design efficiencies*

| $\lambda_0$ | 0 | 1/40 | 5/40 | 10/40 | 10/40 | 10/40 |
|---|---|---|---|---|---|---|
| $\lambda_1$ | 1 | 1 | 1 | 1 | 1/2 | 1/10 |
| $\pi_0$ | 71 | 73 | 77 | 83 | 100 | 100 |
| $\pi_1$ | 97 | 98 | 100 | 100 | 98 | 86 |
| $\pi_2$ | 100 | 100 | 95 | 83 | 80 | 69 |

value of $\lambda_0/\lambda_1$ are in the same interval, then the chosen order is optimal. Next, observe that

$$\sigma_{0\varepsilon}^2 \operatorname{trace}(C_d^L(\pi_q)) = b\left(k - k\lambda_0 - \lambda_1 - \frac{k}{v}(1 - k\lambda_0)\right) + 2b\sum_{p=1}^{q}(\lambda_1\phi^2(p) - \lambda_0).$$

The efficiency of $\pi_q$ may be defined as

$$E_2 = \frac{\operatorname{trace}(C_d^L(\pi_q))}{\operatorname{trace}(C_d^L(\pi_{q^*}))},$$

where $q^*$ is either $s^*$ or $\alpha$, depending on which order is optimal, according to Lemma 3.2. Note that again, $E_2$ does not depend on $b$. It can be shown that the efficiency gets smaller as we move away from the optimal order $\pi_{q^*}$. We will limit our consideration of the magnitude of the efficiencies to a small example.

Let $k = 4$ and $v = 7$. Depending on the value of $\lambda_0/\lambda_1$, an optimal order is either $\pi_0 = \{1, 2, 3, 4\}$, $\pi_1 = \{1, 2, 3, 1\}$ or $\pi_3 = \{1, 2, 2, 1\}$. More precisely, $\pi_0$ is optimal if $\lambda_0/\lambda_1 \geq \phi^2(1) = 9/20$, $\pi_1$ is optimal if $1/20 = \phi^2(2) \leq \lambda_0/\lambda_1 < \phi^2(1) = 9/20$ and $\pi_2$ is optimal if $\lambda_0/\lambda_1 < \phi^2(2) = 1/20$. Table 1 shows the efficiencies (rounded to nearest percentages) of these designs for selected values of $\lambda_0$ and $\lambda_1$.

The conclusion is that we need to be a bit more careful in this case, but that a design that is less extreme (in terms of the value of $q$) is more likely to keep a high efficiency, except possibly for extreme values of $\lambda_0/\lambda_1$.

## 5. Proofs.

PROOF OF LEMMA 3.2.   (i) Suppose $s^* > 0$. For $\pi \in \mathcal{O}$ and $j = 0, 1, \ldots, k$, let us define

(5.1)     $s_j = s_j(\pi) =$ number of symbols that appear $j$ times in $\pi$.

It follows that $v = s_0 + s_1 + \cdots + s_k$, $k = s_1 + 2s_2 + \cdots + ks_k$ and $s = s(\pi) = s_2 + \binom{3}{2}s_3 + \cdots + \binom{k}{2}s_k$. Suppose positions $p_1$, $p_2$ are occupied by symbol $i_1$, positions $p_3$, $p_4$, $p_5$ are occupied by $i_2$ and so on. We can then write

$$
\begin{aligned}
(5.2) \qquad F(\pi) &= -\lambda_0 s - \lambda_1 [\phi(p_1)\phi(p_2) + \phi(p_3)\phi(p_4) \\
&\qquad\qquad\quad + \phi(p_3)\phi(p_5) + \phi(p_4)\phi(p_5) + \cdots] \\
&= -\lambda_0 s - \frac{\lambda_1}{2}[(\phi(p_1) + \phi(p_2))^2 - \phi^2(p_1) - \phi^2(p_2) \\
&\qquad\qquad\quad + (\phi(p_3) + \phi(p_4) + \phi(p_5))^2 \\
&\qquad\qquad\qquad - \phi^2(p_3) - \phi^2(p_4) - \phi^2(p_5) + \cdots] \\
&\leq -\lambda_0 s + \frac{\lambda_1}{2}\left[\sum_{i=1}^{h} \phi^2(p_i)\right],
\end{aligned}
$$

where

$$
(5.3) \qquad h = h(\pi) = 2s_2 + 3s_3 + \cdots + ks_k = k - s_1.
$$

Note that

$$
(5.4) \qquad s - h/2 = \sum_{l=2}^{k} \left(\frac{l(l-1)}{2} - \frac{l}{2}\right)s_l \geq 0.
$$

Hence, we get from (5.2)

$$
(5.5) \qquad F(\pi) \leq -\lambda_0 \frac{h}{2} + \frac{\lambda_1}{2}\left[\sum_{i=1}^{h} \phi^2(p_i)\right] = \frac{1}{2}\sum_{i=1}^{h}[\lambda_1\phi^2(p_i) - \lambda_0].
$$

It follows from (3.9) that an upper bound for $F(\pi)$ is given by

$$
(5.6) \qquad F(\pi) \leq \sum_{p=1}^{s^*}[\lambda_1\phi^2(p) - \lambda_0],
$$

which is attained by the order $\pi_{s^*}$ defined in (3.8).

When $s^* = 0$, $F(\pi_{s^*}) = 0$ and for any other order $\pi$, it follows from (5.5) and (3.9) that $F(\pi) \leq 0$. This establishes part (i) of Lemma 3.2.

(ii) If equality is attained in (5.6) by an order $\hat{\pi}$, then it is clear from (5.4) that $h(\hat{\pi}) = 2s^*$, $s_j(\hat{\pi}) = 0$ for $j \geq 3$, $s_2(\hat{\pi}) = s^*$ and $s_1(\hat{\pi}) = k - 2s^*$. However, for $k > v + s^*$, the number of treatments required by $\hat{\pi}$ is $k - 2s^* + s^* = k - s^* > v$ so that $\hat{\pi}$ does not exist. To prove (ii), let us start by evaluating the possible range of $s(\pi)$ for an optimal order $\pi$. Since $1 < k/v < 2$, it follows from (3.6) that $s(\pi)$ is minimized when $n_i \in \{1, 2\}$; It can hence be shown that $\text{Min}_{\pi \in \mathcal{O}} s(\pi) = \alpha$. On the other hand, if $s(\pi) > k/2$, then the trend-free order $\pi^e$ defined by

$$
\pi^e = \begin{cases} (1, 2, \ldots, k/2, k/2, \ldots, 2, 1), & \text{if } k \text{ is even,} \\ (1, 2, \ldots, (k-1)/2, (k+1)/2, (k-1)/2, \ldots, 2, 1), & \text{if } k \text{ is odd,} \end{cases}
$$

satisfies $s(\pi) > s(\ \pi^e) = \lfloor k/2 \rfloor$ and $T(\pi) \geq T(\pi^e) = -1/2 = \text{Min}_{\pi \in \mathcal{O}}\ T(\pi)$, hence $F(\pi) < F(\pi^e)$. Therefore, without loss of generality, we may assume $\alpha \leq s(\pi) \leq k/2$.

Suppose $\overline{\pi} \in \mathcal{O}$ is arbitrary with $s(\overline{\pi}) = \alpha + g$, where $0 \leq g \leq k/2 - \alpha$. It follows from (5.3) and (5.4) that $s_1(\overline{\pi}) \geq k - 2(\alpha + g)$. Therefore, from (3.7), we get

$$T(\overline{\pi}) = \tfrac{1}{2}\left[\sum_{i=1}^{v} h_i^2(\overline{\pi}) - 1\right] \geq \tfrac{1}{2}\left[\sum_{i:n_i=1} h_i^2(\overline{\pi}) - 1\right] \geq -\sum_{p=1}^{\alpha+g} \phi^2(p).$$

Using (3.3), we obtain $F(\overline{\pi}) \leq \text{Max}_{s(\pi)=\alpha+g}\ F(\pi) = \sum_{p=1}^{\alpha+g}[\lambda_1 \phi^2(p) - \lambda_0]$. Since $(k+1)/2 > \alpha + g \geq \alpha > s^*$, (3.9) implies $\sum_{p=1}^{\alpha+g}[\lambda_1 \phi^2(p) - \lambda_0] \leq \sum_{p=1}^{\alpha}[\lambda_1 \phi^2(p) - \lambda_0]$. It can be shown that $\sum_{p=1}^{\alpha}[\lambda_1 \phi^2(p) - \lambda_0] = F(\pi_\alpha)$, with $\pi_\alpha$ as defined in (3.8). Hence, $F(\overline{\pi}) \leq F(\pi_\alpha)$. This completes the proof of Lemma 3.2.  □

To prove Lemma 3.5, we need a result that is stated and proved below.

LEMMA 5.1.   *Let $k \geq 2v$ and $k$ be even. With $m$ and $t$ as defined in (3.10) and with $u$ denoting an integer, $0 \leq u \leq v$, let $\mathcal{O}_u$ denote the set of all orders with precisely $u$ treatments that have an odd replication. Then, an order $\pi \in \mathcal{O}_u$ with the following properties minimizes $\sum_{i=1}^{v} n_i^2$ over $\mathcal{O}_u$:*

(i) *when $m$ is even,*

$$s_m = v - \frac{u+t}{2}, \qquad s_{m+1} = u, \qquad s_{m+2} = \frac{t-u}{2} \qquad if\ u \leq t,$$

$$s_{m-1} = \frac{u-t}{2}, \qquad s_m = v - u, \qquad s_{m+1} = \frac{u+t}{2} \qquad if\ u > t;$$

(ii) *when $m$ is odd,*

$$s_{m-1} = \frac{v-u-t}{2}, \qquad s_m = u, \qquad s_{m+1} = \frac{v-u+t}{2} \qquad if\ u \leq v - t,$$

$$s_m = \frac{u+v-t}{2}, \qquad s_{m+1} = v - u, \qquad s_{m+2} = \frac{u-v+t}{2} \qquad if\ u > v - t.$$

*Here, the $s_j$'s are the quantities defined in (5.1).*

PROOF.   We will first show that an order that minimizes $\sum_{i=1}^{v} n_i^2$ in $\mathcal{O}_u$ (a "minimizing order") satisfies the following:

(5.7)                    If $s_{j_0} > 0$   and   $s_{j_1} > 0$     then $|j_1 - j_0| \leq 2$.

To see this, suppose $\pi$ is an order such that $s_{j_0} > 0$, $s_{j_1} > 0$ for $j_1 \geq j_0 + 3$. Suppose $n_{i_0} = j_0$ and $n_{i_1} = j_1$. Let $\pi'$ be an order obtained from $\pi$ by only changing two appearances of treatment $i_1$ to treatment $i_0$. For $\pi'$, $n'_{i_1} = j_1 - 2$, $n'_{i_0} = j_0 + 2$, $n'_i = n_i$ for all $i \neq i_0, i_1$, hence $\pi' \in \mathcal{O}_u$. Clearly, $\sum_{i=1}^{v} n'^2_i - \sum_{i=1}^{v} n^2_i = (j_0 + 2)^2 - j_0^2 + (j_1 - 2)^2 - j_1^2 = 4(j_0 - j_1) + 8 \leq -12 + 8 < 0$. Hence, $\pi'$ is "better" than $\pi$. For $\pi'$, $s'_{j_0} = s_{j_0} - 1$, $s'_{j_1} = s_{j_1} - 1$, $s'_{j_0+2} = s_{j_0+2} + 1$, $s'_{j_1-2} = s_{j_1-2} + 1$ and $s'_j = s_j$ for $j \notin \{j_0, j_0 + 2, j_1 - 2, j_1\}$. Repeated application shows that (5.7) must hold for a minimizing order.

Since $k = mv + t$, $0 \leq t < v$, for a minimizing order, we have two possibilities:

$$(5.8) \qquad s_j = 0 \qquad \text{for } j \notin \{m, m+1, m+2\} \text{ or}$$

$$(5.9) \qquad s_j = 0 \qquad \text{for } j \notin \{m-1, m, m+1\}.$$

Suppose $m$ is even, $m \geq 2$. Clearly, $u$ and $t$ are even. For $\pi \in \mathcal{O}_u$, if (5.8) holds, then $s_{m+1} = u$. It follows from

$$(5.10) \quad s_m + s_{m+1} + s_{m+2} = v, \qquad ms_m + (m+1)s_{m+1} + (m+2)s_{m+2} = k$$

that an order $\pi \in \mathcal{O}_u$ with $s_{m+1} = u$ must satisfy

$$(5.11) \qquad s_m = v - \frac{u+t}{2}, \qquad s_{m+1} = u, \qquad s_{m+2} = \frac{t-u}{2}.$$

On the other hand, if (5.9) holds, then $s_{m-1} + s_{m+1} = u$. Identities (5.10) imply that an order $\pi \in \mathcal{O}_u$ with $s_{m-1} + s_{m+1} = u$ must satisfy

$$(5.12) \qquad s_{m-1} = \frac{u-t}{2}, \qquad s_m = v - u, \qquad s_{m+1} = \frac{u+t}{2}.$$

It is clear that when $u < t$, (5.12) cannot hold and when $u > t$, (5.11) cannot hold. When $u = t$, (5.11) and (5.12) both reduce to $s_m = v - t$, $s_{m+1} = t$. This proves Lemma 5.1 for even $m$. The proof for the case of odd $m$, $m \geq 3$, is similar. $\square$

PROOF OF LEMMA 3.5. To prove (i), note that since $|n_i - n_{i'}| \leq 1$, it follows from (3.6) that $\pi^A_{TF}$ minimizes $s(\pi)$ and since $h_i = 0$ for $i = 1, \ldots, v$, it follows from (3.7) that $\pi^A_{TF}$ minimizes $T(\pi)$. Hence, $\pi^A_{TF}$ maximizes $F(\pi)$. The proof of (ii) is similar.

To prove (iii), first consider an order $\pi \in \mathcal{O}_u$. From (3.13), it follows that

$$T(\pi) \geq \frac{1}{2}\left[u\phi^2\left(\frac{k}{2}\right) - 1\right].$$

Case 1: $m$ is even, $m \geq 2$. It follows from Lemma 5.1 that when $\pi \in \mathcal{O}_u$, for $u \leq t$,

$$s(\pi) \geq \frac{1}{2}\left[\left(v - \frac{u+t}{2}\right)m^2 + u(m+1)^2 + \frac{t-u}{2}(m+2)^2 - mv - t\right] = s^{\min} + \frac{t-u}{2},$$

with $s^{\min}$ as defined in (4.2). Similarly, for $u > t$, $s(\pi) \geq s^{\min} + (u - t)/2$. Hence, from (3.3), we get, for $\pi \in \mathcal{O}_u$,

$$(5.13) \qquad F(\pi) \leq -\lambda_0 \left( s^{\min} + \frac{|t - u|}{2} \right) - \frac{\lambda_1}{2} \left[ u\phi^2 \left( \frac{k}{2} \right) - 1 \right].$$

If we denote the upper bound in (5.13) by $F^*(u)$, then $F^*(u) < F^*(t)$ for $u > t$. This implies that $\text{Max}_{0 \leq u \leq v} F^*(u)$ is attained at some $u \leq t$. When $u \leq t$,

$$F^*(u) = -\lambda_0 s^{\min} - \frac{\lambda_0 t}{2} + \frac{\lambda_1}{2} + \frac{u}{2} \left[ \lambda_0 - \lambda_1 \phi^2 \left( \frac{k}{2} \right) \right]$$

$$\leq \begin{cases} F^*(0), & \text{when } \lambda_0 - \lambda_1 \phi^2 \left( \frac{k}{2} \right) < 0, \\ F^*(t), & \text{when } \lambda_0 - \lambda_1 \phi^2 \left( \frac{k}{2} \right) \geq 0. \end{cases}$$

The lemma follows since $F(\pi_{TF}^C) = F^*(0)$ and $F(\pi_{NTF}) = F^*(t)$.

*Case* 2: $m$ is odd, $m \geq 3$. The proof is similar. It can be shown that when $\pi \in \mathcal{O}_u$,

$$F(\pi) \leq -\lambda_0 \left( s^{\min} + \frac{|v - t - u|}{2} \right) - \frac{\lambda_1}{2} \left[ u\phi^2 \left( \frac{k}{2} \right) - 1 \right] = F^{**}(u), \qquad \text{say}.$$

$\text{Max}_{0 \leq u \leq v} F^{**}(u)$ is attained at some $u \leq v - t$. When $u \leq v - t$,

$$F^{**}(u) \leq \begin{cases} F^{**}(0), & \text{when } \lambda_0 - \lambda_1 \phi^2 \left( \frac{k}{2} \right) < 0, \\ F^{**}(v - t), & \text{when } \lambda_0 - \lambda_1 \phi^2 \left( \frac{k}{2} \right) \geq 0. \end{cases}$$

Finally, in this case, $F(\pi_{TF}^C) = F^{**}(0)$ and $F(\pi_{NTF}) = F^{**}(v - t)$. $\quad \square$

DEPARTMENT OF MATHEMATICS, STATISTICS
   AND COMPUTER SCIENCE
UNIVERSITY OF ILLINOIS AT CHICAGO
CHICAGO, ILLINOIS 60607-7045
USA
E-MAIL: dibyen@uic.edu

DEPARTMENT OF STATISTICS
UNIVERSITY OF GEORGIA
ATHENS, GEORGIA 30602
USA
E-MAIL: jstufken@uga.edu